\documentclass[a4paper]{amsart}
\usepackage{graphicx}
\usepackage{amsmath}
\usepackage{latexsym}
\usepackage{amssymb}
\usepackage[all]{xy}
\xyoption{matrix}
\xyoption{arrow}
\begin{document}

\renewcommand{\th}{\operatorname{th}\nolimits}
\newcommand{\rej}{\operatorname{rej}\nolimits}
\newcommand{\extto}{\xrightarrow}
\renewcommand{\mod}{\operatorname{mod}\nolimits}
\newcommand{\ul}{\underline}
\newcommand{\Sub}{\operatorname{Sub}\nolimits}
\newcommand{\ind}{\operatorname{ind}\nolimits}
\newcommand{\Fac}{\operatorname{Fac}\nolimits}
\newcommand{\add}{\operatorname{add}\nolimits}
\newcommand{\soc}{\operatorname{soc}\nolimits}
\newcommand{\Hom}{\operatorname{Hom}\nolimits}
\newcommand{\Rad}{\operatorname{Rad}\nolimits}
\newcommand{\RHom}{\operatorname{RHom}\nolimits}
\newcommand{\uHom}{\operatorname{\underline{Hom}}\nolimits}
\newcommand{\End}{\operatorname{End}\nolimits}
\renewcommand{\Im}{\operatorname{Im}\nolimits}
\newcommand{\Ker}{\operatorname{Ker}\nolimits}
\newcommand{\Coker}{\operatorname{Coker}\nolimits}
\newcommand{\Ext}{\operatorname{Ext}\nolimits}
\newcommand{\op}{{\operatorname{op}}}
\newcommand{\Ab}{\operatorname{Ab}\nolimits}
\newcommand{\id}{\operatorname{id}\nolimits}
\newcommand{\pd}{\operatorname{pd}\nolimits}
\newcommand{\A}{\operatorname{\mathcal A}\nolimits}
\newcommand{\C}{\operatorname{\mathcal C}\nolimits}
\newcommand{\D}{\operatorname{\mathcal D}\nolimits}
\newcommand{\X}{\operatorname{\mathcal X}\nolimits}
\newcommand{\Y}{\operatorname{\mathcal Y}\nolimits}
\newcommand{\F}{\operatorname{\mathcal F}\nolimits}
\newcommand{\Z}{\operatorname{\mathbb Z}\nolimits}
\renewcommand{\P}{\operatorname{\mathcal P}\nolimits}
\newcommand{\T}{\operatorname{\mathcal T}\nolimits}
\newcommand{\G}{\Gamma}
\renewcommand{\L}{\Lambda}
\newcommand{\bdot}{\scriptscriptstyle\bullet}
\renewcommand{\r}{\operatorname{\underline{r}}\nolimits}
\newtheorem{lemma}{Lemma}[section]
\newtheorem{prop}[lemma]{Proposition}
\newtheorem{cor}[lemma]{Corollary}
\newtheorem{theorem}[lemma]{Theorem}
\newtheorem{remark}[lemma]{Remark}
\newtheorem{definition}[lemma]{Definition}

%%%%%%%%%%%%% karins additions:
\newtheorem{ex}[lemma]{Example}
\newtheorem{conj}{Conjecture}

\title[A geometric description of $m$-cluster categories]{A geometric description of $m$-cluster categories}

\author[Baur]{Karin Baur}
\address{Department of Mathematics \\
University of Leicester \\
University Road \\
Leicester LE1 7RH \\
England
}
\email{k.baur@mcs.le.ac.uk}

\author[Marsh]{Robert J. Marsh}
\address{Department of Mathematics \\
University of Leicester \\
University Road \\
Leicester LE1 7RH \\
England
}
\email{rjm25@mcs.le.ac.uk}


\begin{abstract}
We show that the $m$-cluster category of type $A_{n-1}$ is equivalent to a
certain geometrically-defined category of diagonals of a regular $nm+2$-gon.
This generalises a result of Caldero, Chapoton and Schiffler for $m=1$.
The approach uses the theory of translation quivers and their corresponding
mesh categories. We also introduce the notion of the $m$th power of a
translation quiver and show how it can be used to realise the
$m$-cluster category in terms of the cluster category.
\end{abstract}

\maketitle

\section*{Introduction}
Let $n,m\in\mathbb{N}$ and let $\Pi$ be a regular $nm+2$-sided polygon.
We show that a category $\mathcal{C}^m_{A_{n-1}}$ of diagonals can be
associated to $\Pi$ in a natural
way. The objects of $\mathcal{C}^m_{A_{n-1}}$ are the diagonals in
$\Pi$ which divide $\Pi$ into
two polygons whose numbers of sides are congruent to $2$ modulo $m$.
A quiver $\Gamma^m_{A_{n-1}}$ can be defined on the set of such diagonals,
with arrows given by a simple geometrical rule. It is shown that this quiver
is a stable translation quiver in the sense of Riedtmann~\cite{riedtmann} with
translation $\tau$ given by a certain rotation of the polygon.
For a field $k$, the category $\mathcal{C}^m_{A_{n-1}}$ is defined as the
mesh category associated to $(\Gamma^m_{A_{n-1}},\tau)$.

Let $Q$ be a Dynkin quiver of type $A_{n-1}$, and let $D^b(kQ)$ denote
the bounded derived category of finite dimensional $kQ$-modules. Let
$\tau$ denote the Auslander-Reiten translate of $D^b(kQ)$, and let
$S$ denote the shift. These are both autoequivalences of $D^b(kQ)$.
Our main result is that
$\mathcal{C}^m_{A_{n-1}}$ is equivalent to the quotient of $D^b(kQ)$ by
the autoequivalence $\tau^{-1}S^m$. We thus obtain a geometric description
of this category in terms of $\Pi$.

The $m$-cluster category $D^b(kQ)/{\tau^{-1}S^m}$ associated to $kQ$
was introduced in~\cite{keller} and has also been studied
in~\cite{thomas} and~\cite{zhu}.
It is a generalisation of the cluster category defined in~\cite{ccs1}
(for type $A$) and~\cite{bmrrt} (general hereditary case).
Thomas and Zhu show that the $m$-cluster category possesses properties
similar to those of the cluster category (his result applies more generally
to any simply-laced Dynkin quiver).  Keller has shown that it is Calabi-Yau
of dimension $m+1$~\cite{keller}. We remark that such Calabi-Yau categories
have also been studied in~\cite{kellerreiten}.

Our definition is motivated by and is a generalisation of the construction
of the cluster category in type $A$ given in~\cite{ccs1},
where a category of diagonals of a polygon is introduced.
The authors show that this category is equivalent to the cluster category
associated to $kQ$. This can be regarded as the
case $m=1$ here. The aim of the current paper is to generalise the
construction of~\cite{ccs1} to the diagonals arising in the $m$-divisible
polygon dissections considered in~\cite{tzanaki} (see Section~\ref{notation}
for more details). We also remark that a connection between the $m$-cluster
category associated to $kQ$ and the diagonals considered here was given in
\cite{thomas}.

We further show that if $(\Gamma,\tau)$ is any stable translation quiver
then the quiver $\Gamma^m$ with the same vertices but with arrows given by
sectional paths in $\Gamma$ of length $m$ is again a stable translation
quiver with translation given by $\tau^m$.
If $(\Gamma,\tau)$ is taken to be the Auslander-Reiten quiver of the
cluster category of a Dynkin quiver of type $A_{nm-1}$, we show that
$\Gamma^m$ contains $\Gamma^m_{A_{n-1}}$ as a connected component; it follows
that the $m$-cluster category is a full subcategory of the additive
category generated by the mesh category of $\Gamma^m$.

Since $\Gamma$ is known to have a geometric construction, our definition
provides a geometric construction for the additive category generated by
the mesh category of any connected component of $\Gamma^m$. We give an
example to show that this provides a geometric construction for quotients
of $D^b(kQ)$ other than the $m$-cluster category.

\section{Notation and Definitions} \label{notation}

In~\cite{tzanaki}, E. Tzanaki studied an abstract simplicial 
complex obtained by dividing a polygon into smaller polygons. 

We recall the definition of an abstract simplicial complex. 
Let $X$ be a finite set and 
$\triangle\subseteq {\mathcal P}(X)$ a collection of subsets. 
Assume that $\triangle$ is closed under taking subsets (i.e. if 
$A\in \triangle$ and $B\subseteq A$ then $B\in \triangle$). 
Then $\triangle$ is 
an {\em abstract simplicial complex} on the ground set 
$X$. The vertices of $S$ are the single element subsets 
of $\triangle$ (i.e. $\{A\}\in \triangle$). 
The faces 
are the elements of $\triangle$, the facets are the maximal 
among those (i.e. the $A\in \triangle$ such that if 
$A\subseteq B$ 
and $B\in \triangle$ then $A=B$). The dimension of a face $A$ is 
equal to $|A|-1$ (where $|A|$ is the cardinality of $A$). 
The complex is said to be \emph{pure} of dimension $d$ if all its 
facets have dimension $d$. 

Let $\Pi$ be an $nm+2$-gon, $m,n\in\mathbb{N}$, with vertices
numbered clockwise from $1$ to $nm+2$. We regard all operations on
vertices of $\Pi$ modulo $nm+2$. A diagonal $D$ is
denoted by the pair $(i,j)$ (or simply by the pair $ij$ if
$1\leq i,j\leq 9$). Thus $(i,j)$ is the same as $(j,i)$.
We call a diagonal $D$ in $\Pi$ an {\itshape $m$-diagonal} if 
$D$ divides $\Pi$ into an $(mj+2)$-gon and an $(m(n-j)+2)$-gon
where $j=1,\dots,$ $\lceil\frac{n-1}{2}\rceil$. 
Then Tzanaki defines the abstract simplicial complex 
$\triangle=\triangle_{A_{n-1}}^m$ on the $m$-diagonals of 
$\Pi$ as follows. 

The vertices of $\triangle$ are the $m$-diagonals. 
The faces of $\triangle_{A_{n-1}}^m$ are the sets of 
$m$-diagonals which pairwise don't cross. They are 
called {\em $m$-divisible dissections (of $\Pi$)}. Then the 
facets are the maximal collections of such $m$-diagonals. 
Each facet contains exactly $n-1$ elements, so the 
complex $\triangle_{A_{n-1}}^m$ is pure of dimension 
$n-2$. 

The case $m=1$ is the complex whose facets are 
triangulations of an $n+2$-gon. 

\section{A stable translation quiver of diagonals}
\label{se:quiver}

To $\triangle=\triangle_{A_{n-1}}^m$ we associate a category 
along the lines of~\cite{ccs1}. As a first step, we 
associate to the simplicial complex a quiver, called 
$\Gamma_{A_{n-1}}^m$. The vertices of the quiver 
are the $m$-diagonals in the defining polygon $\Pi$, 
i.e. the vertices of $\triangle_{A_{n-1}}^m$. 

The arrows of $\Gamma_{A_{n-1}}^m$ are obtained in 
the following way: 

Let $D$, $D'$ be $m$-diagonals with a common vertex $i$ 
of $\Pi$. Let $j$ and $j'$ be the other endpoints of 
$D$, respectively $D'$. The points $i,j,j'$ divide the 
boundary of the polygon $\Pi$ into three arcs,
linking $i$ to $j$, $j$ to $j'$ and $j'$ to $i$.
(We usually refer to a part of the boundary connecting one
vertex to another as an arc.)
If $D$, $D'$ and the arc from $j$ to $j'$ form an $m+2$-gon 
in $\Pi$ and if furthermore, $D$ can be rotated 
clockwise to $D'$ about the common
endpoint $i$, we draw an arrow from $D$ to $D'$ in $\Gamma^m_{A_{n-1}}$.
(By this we mean that $D$ can be rotated clockwise to the line
through $D'$.) Note that if 
$D$, $D'$ are vertices of the quiver $\Gamma_{A_{n-1}}^m$ 
then there is at most one arrow between them.  

Examples~\ref{n=4m=1} and~\ref{m=2n=4} below illustrate 
this construction. 

We then define an automorphism $\tau_m$ of the quiver: 
let 
$\tau_m:\Gamma_{A_{n-1}}^m\to\Gamma_{A_{n-1}}^m$ 
be the map given by 
$D\mapsto D'$ if $D'$ is obtained from $D$ by 
an anticlockwise rotation through $\frac{2m\pi}{nm+2}$ 
about the centre of the polygon. 
Clearly, $\tau_m$ is a bijective map and a morphism of
quivers.

\begin{figure}[ht] 
\begin{center}
\includegraphics{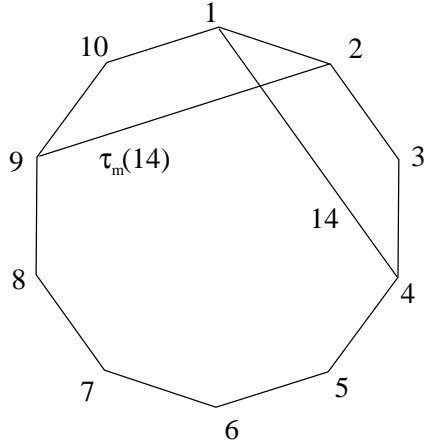}
\caption{The translation $\tau_m$, $\tau_m(14)=92$, where 
$n=4,m=2$}\label{fig1}
\end{center}
\end{figure}

\begin{definition}
(1) A \emph{translation quiver} is a pair $(\Gamma,\tau)$ 
where $\Gamma$ is a locally finite quiver and $\tau:\Gamma'_0\to
\Gamma_0$ is an injective map defined on a subset $\Gamma'_0$ of
the vertices of $\Gamma$ such that for any $X\in\Gamma_0$, $Y\in\Gamma'_0$,
the number of arrows from $X$ to $Y$ is the same as the number of 
arrows from $\tau(Y)$ to $X$. The vertices in $\Gamma_0\setminus \Gamma'_0$
are called \emph{projective}. If $\Gamma'_0=\Gamma_0$ and $\tau$ is
bijective, $(\Gamma,\tau)$ is called a \emph{stable translation quiver}.

(2) A stable translation quiver is said to be \emph{connected} if it is not a 
disjoint union of two non-empty stable subquivers. 
\end{definition}

\begin{prop}\label{prop:stable}
The pair 
$(\Gamma=\Gamma_{A_{n-1}}^m, \tau_m)$ is a stable 
translation quiver. 
\end{prop}

\begin{proof}
By definition, $\tau_m$ is a bijective map from 
$\Gamma$ to $\Gamma$,  and $\Gamma$ is a finite 
quiver. We have to check that the number of  
arrows from $D$ to $D'$ in  
$\Gamma$ is the same as the number 
of arrows from $\tau_m D'$ to $D$. Since there is 
at most one arrow from one vertex to another, 
we only have to see that there is an arrow 
$D\to D'$ if and only if there is an arrow 
$\tau_m D'\to D$.

Assume that there is an arrow $D\to D'$, let $i$ 
be the common vertex of $D$ and $D'$ in the polygon, 
$D=(i,j)$, $D'=(i,j+m)$. Then $\tau_m D'=(i-m,j)$. 
In particular, $j$ is the common vertex of $D$ 
and $\tau_m D'$. Furthermore, we obtain $D$ from
$\tau_m D'$ by a clockwise rotation
about $j$ and these two $m$-diagonals form an $m+2$-gon 
together with an arc from $i-m$ to $i$, hence there is 
an arrow $\tau_m D'\to D$. 

The converse follows with the same reasoning. 
\begin{figure}[ht] 
\begin{center}
\includegraphics{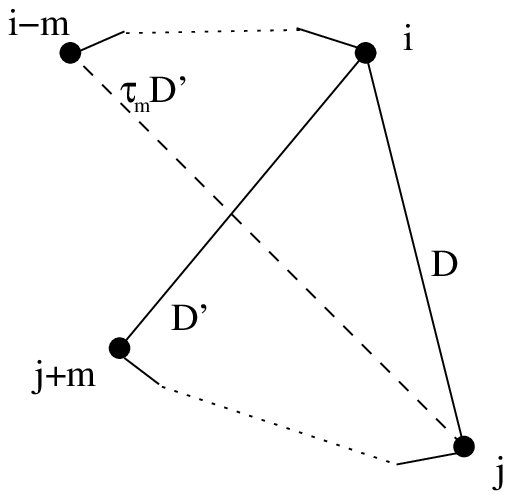}
\caption{$D\to D'$ $\Longleftrightarrow$ $\tau_m D'\to D$}\label{fig2}
\end{center}
\end{figure} 
\end{proof}

\begin{prop}\label{prop:connected}
$(\Gamma,\tau_m)$ is a connected stable translation quiver.
\end{prop}

\begin{proof}
Note that every vertex of $\Pi$ is incident with some element of
any given $\tau_m$-orbit of $m$-diagonals:
any $m$-diagonal is of the form 
$(i,i+km+1)$ and 
\begin{eqnarray*}
\tau_m^{k-n}(i,i+km+1) & = & (i+(n-k)m,i+nm+1) \\ 
  & = & (i+(n-k)m,i-1). 
\end{eqnarray*}

Assume that $\Gamma$ is the disjoint union of two 
non-empty stable subquivers. So there exist $m$-diagonals 
$D=(i,j)$ and $D'=(i',j')$ that cannot be 
connected by any path in $\Gamma$.
Without loss of generality, $j<j'$. 
After rotating $D'$ using $\tau_m$ we can assume 
that $i=i'$. By assumption, $j'\neq j+rm$ for any $r$. 
The diagonal $D$ can be rotated clockwise about $i$ 
to another $m$-diagonal 
$D''=(i,j'')$ such that $j'=j''+s$ with $0<s<m$. 
Since $D''$ is an $m$-diagonal, the arc from $i$ to $j''$ not including $j'$ 
together with $D''$ bounds a $(um+2)$-gon for some $u$. 
But then the arc from $i$ to $j'$ including $j''$ together with the
diagonal $D'$ bound a $(um+2+s)$-gon where $um+2<um+2+s<(u+1)m+2$. 
Hence $D'$ cannot be an $m$-diagonal. 
\end{proof}

In the examples below we draw the quiver 
associated to the complex $\triangle_{A_{n-1}}^m$ 
in the standard way of Auslander-Reiten theory: 
the vertices and arrows are arranged so that 
the translation $\tau_m$ is a shift to the left. 
We indicate it by dotted lines.

\begin{ex}\label{n=4m=1}
Let $n=4$, $m=1$, i.e. $\Pi$ is a $6$-gon. The 
rotation group given by rotation about the centre 
of $\Pi$ through 
$k\times \frac{\pi}{3}$ degrees ($k=1,\dots,5$) 
acts on the facets of $\triangle_{A_3}^1$. 
There are four orbits, 
${\mathcal O}_{\{13,14,15\}}$ of size $6$, 
${\mathcal O}_{\{13,14,46\}}$ and 
${\mathcal O}_{\{13,36,46\}}$ of size $3$ and 
${\mathcal O}_{\{13,15,35\}}$ with two elements,
making a total of $14$ elements.

The vertices of the quiver $\Gamma_{A_3}^1$ are the nine
$1$-diagonals $\{13,14,15,24,25,26,35,36,46\}$ and 
we draw the quiver as follows: 

$$
{\small
\xymatrix@-4mm{
%\xymatrix{
46\ar[rd]\ar@{.}[rr] & & 15\ar[rd]\ar@{.}[rr] & 
& 26\ar[rd]\ar@{.}[rr] & & 13\ar[rd] & & \\
\ar@{.}[r] 
 & 14\ar[rd]\ar[ru]\ar@{.}[rr] &  
& 25\ar[rd]\ar[ru]\ar@{.}[rr] &  
& 36\ar[rd]\ar[ru]\ar@{.}[rr] & 
 & 14\ar@{-->}[rd]\ar@{-->}[ru]\ar@{.}[r] &\\
13\ar[ru]\ar@{.}[rr] & 
 & 24\ar[ru]\ar@{.}[rr] &  
 & 35\ar[ru]\ar@{.}[rr] & & 46\ar[ru] & &
}}
$$

\end{ex}

\begin{ex}\label{m=2n=4}
Let $m=2$ and $n=4$, i.e. $\Pi$ is a $10$-gon. The 
rotation group is generated by the rotation 
about the centre of $\Pi$ through $k\times\frac{2\pi}{5}$ 
degrees ($k=1,\dots,9$) and acts on the facets of 
$\triangle_{A_3}^2$. The orbits are 
${\mathcal O}_{\{14,16,18\}}$,
${\mathcal O}_{\{14,18,47\}}$,
${\mathcal O}_{\{18,38,47\}}$ and
${\mathcal O}_{\{47,38,39\}}$ of size $10$, and
${\mathcal O}_{\{14,16,69\}}$,
${\mathcal O}_{\{14,49,69\}}$ and
${\mathcal O}_{\{29,38,47\}}$ of size 5,
making a total of $55$ elements.
The vertices of $\Gamma_{A_3}^2$ are the fifteen $2$-diagonals 

$\{ 14,16,18,25,27,29,36,38,(3,10),47,49,58,(5,10),69,(7,10)\}$ 
 
\noindent
and the quiver is: 
$$
{\small
\xymatrix@-5mm{
%\xymatrix{
69\ar[rd]\ar@{.}[rr] & & 18\ar[rd]\ar@{.}[rr] & 
& 3,10\ar[rd]\ar@{.}[rr] & & 25\ar[rd]\ar@{.}[rr] & 
 & 47\ar@{.}[rr]\ar[rd] & & 69\ar[rd] & & \\
\ar@{.}[r] & 
 16\ar[rd]\ar[ru]\ar@{.}[rr] &  
& 38\ar[rd]\ar[ru]\ar@{.}[rr] &  
& 5,10\ar[rd]\ar[ru]\ar@{.}[rr] & 
 & 27\ar[rd]\ar[ru]\ar@{.}[rr] & 
 & 49\ar@{.}[rr]\ar[rd]\ar[ru]
 & & 16\ar@{-->}[rd]\ar@{-->}[ru]\ar@{.}[r]& \\
14\ar[ru]\ar@{.}[rr] & 
 & 36\ar[ru]\ar@{.}[rr] &  
 & 58\ar[ru]\ar@{.}[rr] & & 7,10\ar[ru]\ar@{.}[rr] & & 
29\ar@{.}[rr]\ar[ru] & & 14\ar[ru] & & 
}}
$$

\end{ex} 

\section{$m$-Cluster Categories}

Let $G$ be a simply-laced Dynkin diagram with vertices $I$.
Let $Q$ be a quiver with underlying
graph $G$, and let $k$ be an algebraically-closed field. Let $kQ$ be the
corresponding path algebra. Let $D^b(kQ)$ denote the bounded derived
category of finitely generated $kQ$-modules, with shift denoted by $S$,
and Auslander-Reiten translate given by $\tau$. It is known that
$D^b(kQ)$ is triangulated, Krull-Schmidt and has almost-split triangles
(see~\cite{happel}).
Let $\mathbb{Z}Q$ be the stable translation quiver associated to $Q$, with
vertices $(n,i)$ for $n\in\mathbb{Z}$ and $i$ a vertex of $Q$. For every
arrow $\alpha:i\to j$ in $Q$ there are arrows $(n,i)\to
(n,j)$ and $(n,j)\to (n+1,i)$ in $\mathbb{Z}Q$, for all
$n\in\mathbb{Z}$. Together with the translation $\tau$, taking $(n,i)$ to
$(n-1,i)$, $\mathbb{Z}Q$ is a stable translation quiver. We note that
$\mathbb{Z}Q$ is independent of the orientation of $Q$ and can thus be
denoted $\mathbb{Z}G$.

We recall the notion of the mesh category of a stable translation 
quiver  with no multiple arrows (the mesh category is 
defined for 
a general translation quiver but we shall not need that here).
Recall that for a quiver $\Gamma$, $k\langle \Gamma \rangle$ denotes the
path category on $\Gamma$, with morphisms given by arbitrary $k$-linear
combinations of paths.

\begin{definition}
Let $(\Gamma,\tau)$ be a stable translation quiver with no multiple 
arrows. 
Let $Y$ be a vertex of $\Gamma$ and let $X_1,\dots,X_k$ be all the 
vertices with arrows to $Y$, denoted $\alpha_i:X_i\to Y$. Let
$\beta_i:\tau(Y)\to X_i$ be the corresponding arrows from $\tau(Y)$ to
$X_i$ ($i=1,\dots,k$).  
Then the {\itshape mesh ending at $Y$} is defined to be the quiver consisting 
of the vertices $Y,\tau(Y),X_1,\dots,X_k$ and the arrows $\alpha_1,\alpha_2,
\ldots ,\alpha_k$ and $\beta_1,\beta_2,\ldots ,\beta_k$.
The {\it mesh relation at Y} is defined to be  
\[
m_Y:=\sum_{i=1}^k \beta_i\alpha_i\in \Hom_{k\langle \Gamma \rangle}(\tau(Y),Y)
\]
Let $I_m$ be the ideal in $k\langle \Gamma\rangle$ generated by the mesh
relations $m_Y$ where $Y$ runs over all vertices of $\Gamma$. 

Then the {\itshape mesh category} of $\Gamma$ is defined as the 
quotient $k\langle\Gamma\rangle/I_m$.
\end{definition}

For an additive category $\varepsilon$, denote by $\ind\varepsilon$ the
full subcategory of indecomposable objects.
Happel~\cite{happel} has shown that $\ind D^b(kQ)$ is equivalent to the
mesh category of $\mathbb{Z}Q$, from which it follows that it is independent
of the orientation of $Q$. Its Auslander-Reiten quiver is $\mathbb{Z}G$.

For $m\in\mathbb{N}$, we denote by $\mathcal{C}^m_G$ the
$m$-cluster category associated to the Dynkin diagram $G$, so
$$\mathcal{C}^m_G=\frac{D^b(kQ)}{F_m},$$
where $Q$ is any orientation of $G$ and $F_m$ is the autoequivalence
$\tau^{-1}\circ S^m$ of $D^b(kQ)$. This was introduced by
Keller~\cite{keller} and has been studied by
Thomas~\cite{thomas} and Zhu~\cite{zhu}.
It is known that $\mathcal{C}^m_G$ is triangulated~\cite{keller},
Krull-Schmidt
and has almost split triangles~\cite[1.2,1.3]{bmrrt}. Let $\varphi_m$
denote the automorphism of $\mathbb{Z}G$
induced by the autoequivalence $F_m$. The Auslander-Reiten
quiver of $\mathcal{C}^m_G$ is the quotient $\mathbb{Z}G/\varphi_m$, and
$\ind\mathcal{C}^m_G$ is equivalent to the mesh category of
$\mathbb{Z}G/\varphi_m$.

\section{Proof of the geometric description}

Our main aim in this section is to show that, if $G$ is of type $A_{n-1}$,
then $\ind\mathcal{C}^m_G$ is equivalent to the mesh category
$\mathcal{D}^m_{A_{n-1}}$ of the stable translation quiver $\Gamma^m_{A_{n-1}}$
defined in the previous section. From the previous section we can see that it
is enough to show that, as translation quivers, $\mathbb{Z}G/\varphi_m$ is
isomorphic to $\Gamma^m_{A_{n-1}}$.

\subsection{$m$-coloured almost positive roots and $m$-diagonals}

For $\Phi$ a root system, with positive roots $\Phi^+$ and simple
roots $\alpha_1,\alpha_2,\ldots ,\alpha_n$, let $\Phi_{\geq -1}^m$ denote the
set of $m$-coloured almost positive roots (see~\cite{fominreading}).
An element of $\Phi_{\geq -1}^m$ is either a $m$-coloured positive root
$\alpha^k$ where $\alpha\in\Phi^+$ and $k\in\{1,2,\ldots ,m\}$ or a negative
simple root $-\alpha_i$ for some $i$ which we regard as having colour $1$
for convenience (it is thus also denoted $-\alpha_i^1$).
Fomin-Reading~\cite{fominreading}
show that there is a one-to-one correspondence between $m$-diagonals of
the regular $nm+2$-gon $\Pi$ and $\Phi^m_{\geq -1}$ when $\Phi$ is of type
$A_{n-1}$. We now recall this correspondence.

\begin{figure}[h] 
\begin{center}
\includegraphics{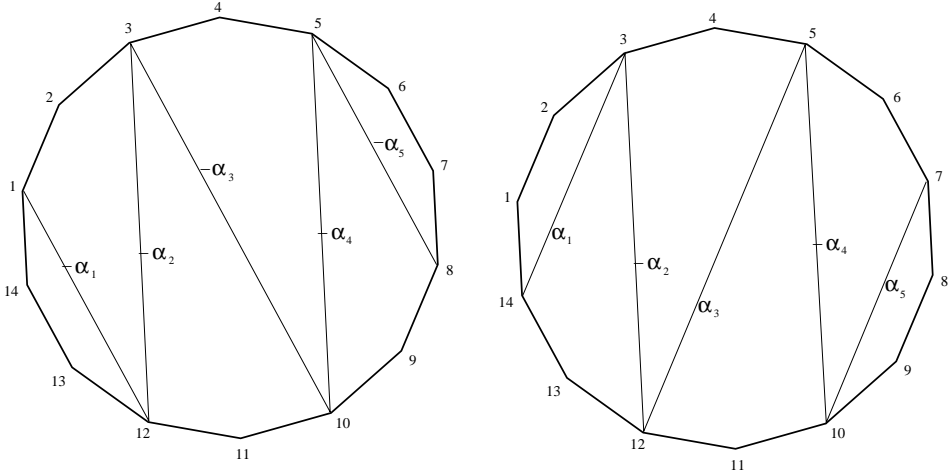}
\caption{$m$-snake and opposite $m$-snake for  
$n=6,m=2$}\label{m-snake}
\end{center}
\end{figure}

Recall that $R_m$ denotes the anticlockwise rotation of $\Pi$ taking vertex $i$
to vertex $i-1$ for $i\geq 2$, and vertex $1$ to vertex $nm+2$.
For $1\leq i\leq \frac{n}{2}$, the negative simple root -$\alpha_{2i-1}$
corresponds to the diagonal $((i-1)m+1,(n-i)m+2)$. For
$1\leq i\leq \frac{n-1}{2}$, the negative simple root $-\alpha_{2i}$
corresponds to the diagonal $(im+1,(n-i)m+2)$.
Together, these diagonals form what is known as the \emph{$m$-snake},
cf. Figure~\ref{m-snake}. 
For $1\leq i\leq j\leq n$, there are exactly $m$ $m$-diagonals intersecting
the diagonals labelled $-\alpha_i,-\alpha_{i+1},\ldots ,-\alpha_j$ and
no other diagonals labelled with negative simple roots. These diagonals
are of the form $D,R_m^1(D),\ldots ,R_m^{m-1}(D)$ for some diagonal $D$,
and $\alpha^k$ corresponds to $R_m^{k-1}(D)$ for $k=1,2,\ldots ,m$, where
$\alpha$ denotes the positive root $\alpha_i+\cdots +\alpha_j$.
For an $m$-coloured almost positive root $\beta^k$,
we denote the corresponding diagonal by $D(\beta^k)$.

It is clear that, for $1\leq i\leq \frac{n}{2}$, the coloured
root $\alpha_{2i-1}^1$ corresponds to the diagonal $(im+1,(n+1-i)m+2)$.
Also, the diagonals $D(-\alpha_i)$,
for $i$ even, together with $D(\alpha_j^1)$, for $j$ odd, form a
`zig-zag' dissection of $\Pi$ which we call the
\emph{opposite $m$-snake}, cf. Figure~\ref{m-snake}.

Let $I=I^+\cup I_-$ be a decomposition of the vertices $I$ of $G$ so that
there are no arrows between vertices in $I^+$ or between vertices in $I^-$;
such a decomposition exists because $G$ is bipartite. For type $A_{n-1}$,
we take $I_+$ to be the even-numbered vertices and $I_-$ to be the odd-numbered
vertices.

Let $R_m:\Phi^m_{\geq -1}\to \Phi^m_{\geq -1}$ be the bijection
introduced by Fomin-Reading~\cite[2.3]{fominreading}. This is defined
using the involutions~\cite{fominzelevinsky2}
$\tau_{\pm}:\Phi_{\geq -1}\to
\Phi_{\geq -1}$ given by
$$\tau_{\varepsilon}(\beta)=\left\{ \begin{array}{cc}
\alpha & \mbox{if\ }\beta=-\alpha_i,\mbox{\ for\ }i\in I_{-\varepsilon} \\
\left(\prod_{i\in I_{\varepsilon}}s_i\right)(\beta) & \mbox{otherwise}.
\end{array}\right.
$$
Then, for $\beta^k\in \Phi^m_{\geq -1}$, we have
$$R_m(\beta^k)=\left\{ \begin{array}{cc}
\beta^{k+1} & \mbox{if\ }\alpha\in\Phi^+\mbox{\ and\ }k<m, \\
((\tau_-\tau_+)(\beta))^1 & \mbox{otherwise}.
\end{array}\right.$$

\begin{lemma} [FOMIN-READING] \label{frcompatible}
For all $\beta_k\in\Phi^m_{\geq -1}$, we have:
$D(R_m(\beta_k))=R_mD(\beta_k)$.
\end{lemma}

\begin{proof}
See the discussion in~\cite[4.1]{fominreading}.
\end{proof}

\subsection{Indecomposable objects in the $m$-cluster category and
$m$-diagonals}

Let $Q_{alt}$ denote the orientation of $G$ obtained by orienting
every arrow to go from a vertex in $I_+$ to a vertex in $I_-$, so that the
vertices in $I_+$ are sources and the vertices in $I_-$ are sinks.

For a positive root $\alpha$, let $V(\alpha)$ denote the corresponding
$kQ_{alt}$-module, regarded as an indecomposable object in $D^b(kQ_{alt})$.
Then it is clear from the definition that the indecomposable objects in
$\mathcal{C}^m_G$ are the objects $S^{k-1}V(\beta)$ for $k=1,2,\ldots ,m$
and $\alpha\in\Phi^+$ and $S^{-1}I_i$ for $I_i$ an indecomposable injective
$kQ_{alt}$-module corresponding to the vertex $i\in I$ (all regarded
as objects in the $m$-cluster category). Following Thomas~\cite{thomas}
or Zhu~\cite{zhu},
we define $V(\alpha^k)$ to be $S^{k-1}V(\beta)$ for $k=1,2,\ldots ,m$,
$\alpha\in\Phi^+$, and $V(-\alpha_i)=S^{-1}I_i$ for $i\in I$.

We have:
\begin{lemma} [THOMAS,ZHU] \label{thomaszhucompatible}
For all $\beta_k\in\Phi_{\geq -1}^m$, $V(R_m\beta^k)\cong SV(\beta_k)$,
where $S$ denotes the autoequivalence of $\mathcal{C}^m_G$ induced by the
shift on $D^b(kQ)$.
\end{lemma}

\begin{proof}
See~\cite[Lemma 2]{thomas} or~\cite[3.8]{zhu}.
\end{proof}

\section{An isomorphism of stable translation quivers}

From the previous two sections, we see that in type $A_{n-1}$, we
have a bijection $D$ from $\Phi^m_{\geq -1}$ to the set of $m$-diagonals
of $\Pi$ and a bijection $V$ from $\Phi^m_{\geq -1}$ to the objects of
$\ind\mathcal{C}^m_{A_{n-1}}$ up to isomorphism,
i.e. to the vertices of the Auslander-Reiten quiver of
$\mathcal{C}^m_{A_{n-1}}$.
Composing the inverse of $D$ with $V$ we obtain a bijection $\psi$ from
the set of $m$-diagonals of $\Pi$ to $\ind\mathcal{C}^m_{A_{n-1}}$.

\begin{lemma} \label{psicompatible}
For every $m$-diagonal $D$ of $\Pi$, we have that
$$\psi(R_m(D))\cong S\psi(D),$$ and therefore that
$$\psi(\tau_m(D))\cong \tau_m(\psi(D)).$$
\end{lemma}

\begin{proof}
The first statement follows immediately from Lemmas~\ref{frcompatible}
and~\ref{thomaszhucompatible}. We can deduce from this that
$\psi(\tau_m(D))=\psi(R_m^m(D))=S^m\psi(D)$ and thus obtain the second
statement, since $S^m$ coincides with $\tau_m$ on every
indecomposable object of $\mathcal{C}^m_{A_{n-1}}$ by the definition
of this category.
\end{proof}

It remains to show that $\psi$ and $\psi^{-1}$ are morphisms of quivers.

\begin{lemma}\label{lm:m-arrows}
\begin{itemize}
\item For $1\leq i\leq \frac{n-1}{2}$, there is an arrow in
$\Gamma^m_{A_{n-1}}$ from $D(-\alpha_{2i-1})$ to $D(-\alpha_{2i})$.
\item For $1\leq i\leq \frac{n-1}{2}$, there is an arrow in
$\Gamma^m_{A_{n-1}}$ from $D(-\alpha_{2i+1})$ to $D(-\alpha_{2i})$.
\item For $1\leq i\leq \frac{n}{2}$, there is an arrow in
$\Gamma^m_{A_{n-1}}$ from $D(-\alpha_{2i})$ to $D(\alpha_{2i-1}^1)$.
\item For $1\leq i\leq \frac{n-2}{2}$, there is an arrow in
$\Gamma^m_{A_{n-1}}$ from $D(-\alpha_{2i})$ to $D(\alpha_{2i+1}^1)$.
\end{itemize}
These are the only arrows amongst the diagonals
$D(-\alpha_i)$ and $D(\alpha_j^1)$, for $1\leq i,j\leq n-1$, with
$j$ odd, in $\Gamma^m_{A_{n-1}}$.
\end{lemma}

\begin{proof}
We firstly note that, for $1\leq i\leq \frac{n-1}{2}$, the diagonals
corresponding to the negative simple roots $-\alpha_{2i-1}$ and
$-\alpha_{2i}$, together with an arc of the boundary containing vertices
$(i-1)m+1,\ldots ,im+1$, bound an $m+2$-gon. The other vertex is
numbered $(n-i)m+2$. Furthermore, $D(-\alpha_{2i-1})$ can be rotated clockwise
about the common end point $(n-i)m+2$ to $D(-\alpha_{2i})$, so there is an
arrow in $\Gamma^m_{A_{n-1}}$ from $D(-\alpha_{2i-1})$ to $D(-\alpha_{2i})$.

Similarly, for $1\leq i\leq \frac{n-2}{2}$, the diagonals
corresponding to the negative simple roots $-\alpha_{2i}$ and $-\alpha_{2i+1}$,
together with an arc of the boundary containing vertices
$(n-i-1)m+2,\ldots ,(n-i)m+2$, bound an $m+2$-gon (with the other vertex being
numbered $im+1$), and $D(-\alpha_{2i+1})$ can be rotated clockwise
about the common end point $im+1$ to $D(-\alpha_{2i})$, so there is an
arrow in $\Gamma^m_{A_{n-1}}$ from $D(-\alpha_{2i+1})$ to $D(-\alpha_{2i})$.

We have observed that, for $1\leq i\leq \frac{n}{2}$, the coloured
root $\alpha_{2i-1}^1$ corresponds to the diagonal $(im+1,(n+1-i)m+2)$.
Consideration of the $m+2$-gon with
vertices $(n-i)m+2,\ldots ,(n+1-i)m+2$ and $im+1$ shows that there is an
arrow from $D(-\alpha_{2i})$ to $D(\alpha_{2i-1}^1)$. For $1\leq i\leq 
\frac{n-2}{2}$, consideration
of the $m+2$-gon with vertices $im+1,\ldots (i+1)m+1$ and $(n-i)m+2$
shows that there is an arrow from $D(-\alpha_{2i})$ to $D(\alpha_{2i+1}^1)$.

The statement that these are the only arrows amongst the diagonals
considered is clear.
\end{proof}

The following follows from the well-known structure of the
Auslander-Reiten quiver of $D^b(kQ)$.

\begin{lemma}\label{lm:A-R-arrows}
\begin{itemize}
\item For $1\leq i\leq \frac{n-1}{2}$, there is an arrow in the
Auslander-Reiten quiver of $\mathcal{C}^m_{A_{n-1}}$ from
$I_{2i-1}[-1]$ to $I_{2i}[-1]$.
\item For $1\leq i\leq \frac{n-2}{2}$,
there is an arrow from $I_{2i+1}[-1]$ to $I_{2i}[-1]$.
\item For $1\leq i\leq \frac{n}{2}$, there is an arrow from $I_{2i}[-1]$ to
$P_{2i-1}$.
\item For $1\leq i\leq \frac{n-2}{2}$ there is an arrow from
$I_{2i}[-1]$ to $P_{2i+1}$.
\end{itemize}
These are the only arrows amongst the vertices $I_i[-1]$ and $P_j$ for
$1\leq i,j\leq n-1$, with $j$ odd, in the Auslander-Reiten quiver of
$\mathcal{C}^m_{A_{n-1}}$.
\end{lemma}

\begin{prop} \label{graphiso}
The map $\psi$ from $m$-diagonals in $\Pi$ to indecomposable objects in
$\mathcal{C}^m_{A_{n-1}}$ is an isomorphism of quivers.
\end{prop}

\begin{proof}
Suppose that $D,E$ are $m$-diagonals in $\Pi$ and that there is an arrow
from $D$ to $E$. Write $D=D(\beta^k)$ and $E=D(\gamma^l)$ for coloured
roots $\beta^k$ and $\gamma^l$. Then $V:=\psi(D)=V(\beta^k)$ and
$W:=\psi(E)=V(\gamma^l)$ are corresponding vertices in the Auslander-Reiten
quiver of $\mathcal{C}^m_{A_{n-1}}$. Since there is an arrow from $D$ to $E$, there
is an $m+2$-gon bounded by $D$ and $E$ and an arc of the boundary of $\Pi$.

Since $D$ is an $m$-diagonal, on the side of $D$ not in the $m+2$-gon,
there is an $dm+2$-gon bounded by $D$ and an arc of the boundary of $\Pi$
for some $d\geq 1$.
Similarly, since $E$ is an $m$-diagonal, on the side of $E$ not in the
$m+2$-gon, there is an $em+2$-gon bounded by $D$ and an arc of the boundary
of $\Pi$, for some $e\geq 1$. It is clear that each of these polygons can be
dissected by an $m$-snake such that, together with $D$ and $E$, we obtain
a `zig-zag' dissection $\chi$ of $\Pi$. Let $v$ be one of its endpoints.
The other endpoint must be $v-m-1$ or $v+m+1$ (modulo $nm+2$).

In the first case, we have that for some $t\in\mathbb{Z}$, $R_m^t(v)=1$
and $R_m^t$ applied to $\chi$ is the $m$-snake. In the second case,
we have that, for some $t\in\mathbb{Z}$, $R_m^t(v)=nm+2$ and $R_m^t$
applied to $\chi$ is the opposite $m$-snake.
It follows from Lemma~\ref{lm:A-R-arrows}
that there is an arrow from
$R_m^t(V)$ to $R_m^t(W)$ in the Auslander-Reiten quiver of $\mathcal{C}^m_{A_{n-1}}$,
and hence from $V$ to $W$.

Conversely, suppose that $V,W$ are vertices of the Auslander-Reiten
quiver of $\mathcal{C}^m_{A_{n-1}}$ and that there is an arrow from $V$ to $W$.
We can write $V=V(\beta^k)$ and $W=V(\gamma^l)$ for coloured roots
$\beta^k$ and $\gamma^l$. Let $D:=\psi^{-1}(V)=D(\beta^k)$ and let
$E:=\psi^{-1}(W)=D(\gamma^l)$.
It is clear that $\tau_m^u(V)\cong I_i[-1]$ for some $i$.
By Lemma~\ref{lm:A-R-arrows}, we must
have that either $\tau_m^u(W)\cong I_{i\pm 1}[-1]$ or $\tau_m^u(W)\cong
P_{i\pm 1}$. In the latter case we must have that $i$ is even.
Note that $S^{um}(V)\cong \tau_m^u(V)$ and $S^{um}(W)\cong \tau_m^u(W)$. 
It follows from Lemmas~\ref{psicompatible} and~\ref{lm:m-arrows} that there is 
an arrow from $R_m^{um}(D)$ to
$R_m^{um}(E)$ in $\Gamma^m_{A_{n-1}}$, and thus from $D$ to $E$.

It follows that $\psi$ is an isomorphism of quivers.
\end{proof}

\begin{prop}
There is an isomorphism $\psi$ of translation quivers between the stable
translation quiver $\Gamma^m_{A_{n-1}}$ of $m$-diagonals and the
Auslander-Reiten quiver of the $m$-cluster category $\mathcal{C}^m_{A_{n-1}}$.
\end{prop}

\begin{proof}
This now follows immediately from Proposition~\ref{graphiso} and
Lemma~\ref{psicompatible}.
\end{proof}

We therefore have our main result.

\begin{theorem}\label{thm:m-cluster}
The $m$-cluster category $\mathcal{C}^m_{A_{n-1}}$ is equivalent to the
additive category generated by the mesh category of the stable translation
quiver $\Gamma^m_{A_{n-1}}$ of $m$-diagonals.
\end{theorem}

We remark that a connection between the $m$-cluster category and the
$m$-diagonals has been given in~\cite{thomas}. In particular, Thomas gives
an interpetation of $\Ext$-groups in the $m$-cluster category in terms of
crossings of diagonals. However, Thomas does not give a construction of
the $m$-cluster category using diagonals.

\section{The $m$-th power of a translation quiver} \label{se:mpower}

In this section we define a new category in 
natural way in which the $m$-cluster category 
$\mathcal{C}^m_{A_{n-1}}$ will appear as a full subcategory. 
We start with a translation quiver $\Gamma$ 
and define its $m$-th power. 

Let $\Gamma$ be a translation quiver with 
translation $\tau$. 

Let $\Gamma^m$ be the quiver whose objects 
are the same as the objects of $\Gamma$ and 
whose arrows are the sectional paths of 
length $m$. 
A path $(x=x_0\to x_1\to\dots\to x_{m-1}\to x_m=y)$ in $\Gamma$
is said to be {\it sectional}
if $\tau x_{i+1}\neq x_{i-1}$ for $i=1,\dots,m-1$ (for which $\tau x_{i+1}$
is defined) (cf.~\cite{ringel}). 
Let $\tau^m$ be the $m$-th power 
of the translation, i.e. 
$\tau^m=\tau\circ\tau\circ\dots\circ\tau$ ($m$ times). Note that
the domain of defininiton of $\tau^m$ is a subset of the domain
of definition $\Gamma'_0$ of $\tau$.

Recall that a translation quiver is said to be \emph{hereditary}
(see~\cite{ringel}) if:
\begin{itemize}
\item for any non-projective vertex $z$, there is an arrow
from some vertex $z'$ to $z$;
\item there is no (oriented) cyclic path of length at least one
containing projective vertices, and
\item If $y$ is a projective vertex and there is an arrow $x\to y$,
then $x$ is projective.
\end{itemize}
The last condition is what we need to ensure that $(\Gamma^m,\tau^m)$
is again a translation quiver:

\begin{theorem} \label{mpower}
Let $(\Gamma,\tau)$ be a translation quiver such that
if $y$ is a projective vertex and there is an arrow $x\to y$,
then $x$ is projective. Then $(\Gamma^m,\tau^m)$ is a translation quiver.
\end{theorem}

\begin{proof}
We prove the following statement by induction on $m$:

\noindent
Suppose that there is a sectional path 
\[
x=x_0\to x_1\to\dots\to x_m=y
\]
in $\Gamma$ and $\tau^my$ is defined. Then $\tau^i x_i$ is defined
for $i=0,1,\ldots ,m$ and there is a sectional path
\[
\tau^my=\tau^mx_m\to\tau^{m-1}x_{m-1}\to\dots
\to\tau x_1\to x=x_0
\]
in $\Gamma$.
Furthermore, if the multiplicities of arrows between consecutive vertices
in the first path are $k_1,k_2,\ldots ,k_m$, the multiplicities of
arrows between consecutive vertices in the second path are
$k_m,k_{m-1},\ldots ,k_1$.

This is clearly true for $m=1$, since $\Gamma$ is a translation
quiver. Suppose it is true for $m-1$, and that
$$x=x_0\to x_1\to \cdots \to x_m=y$$
is a sectional path in $\Gamma$. Since $\tau^{m-1}x_m$ is defined, we
can apply induction to the section path:
$$x_1\to x_2\to \cdots \to x_m$$
to obtain that $\tau^{i-1} x_i$ is defined for $i=1,2,\ldots ,m$ and that
there is a sectional path
$$\tau^{m-1}x_m\to \tau^{m-2}x_{m-1}\to \cdots x_1$$
in $\Gamma$, with multiplicities $k_2,k_3,\ldots ,k_m$.
As $\tau^m x_m$ is defined, $\tau^{m-1}x_m$ is not projective, and it
follows that $\tau^{i-1}x_i$ is not projective for $i=1,2,\ldots ,m$
by our assumption. Therefore $\tau^i x_i$ is defined for $i=1,2,\ldots ,m$.
For $i=2,3,\ldots ,m$, there are $k_i$ arrows from $\tau^{i-1}x_i$ to
$\tau^{i-2}x_{i-1}$. Therefore there are $k_i$ arrows from
$\tau^{i-1}x_{i-1}$ to $\tau^{i-1}x_i$. Thus there are $k_i$ arrows
from $\tau^i x_i$ to $\tau^{i-1}x_{i-1}$. As there are $k_1$ arrows
from $x_0$ to $x_1$, there are $k_1$ arrows from $\tau x_1$ to $x_0$.
If $\tau (\tau^i x_i)=\tau^{i+2}x_{i+2}$ for some $i$ then
$x_i=\tau x_{i+2}$, contradicting the fact that $x_0\to x_1\to
\cdots \to y$ is sectional. It follows that
$$\tau^m x_m\to \tau^{m-1}x_{m-1}\to \cdots \to x_0=x$$
is a sectional path with multiplicities of arrows $k_1,k_2,\ldots ,k_m$
as required.

It follows that the number of sectional paths with sequence of vertices
$x_0,x_1,\ldots ,x_m$ is less than or equal to the number of sectional
paths with sequence of vertices
$\tau^m y=\tau^mx_m,\tau^{m-1}x_{m-1},\ldots ,\tau x_1,x_0=x$.

Suppose that
\[
x=x'_0\to x'_1\to  \cdots \to x'_m=y
\]
is a sectional path from $x$ to $y$ with a different sequence of vertices.
Then $x_i\not=x'_i$ for some $i$, $0<i<m$. It follows that
$\tau^ix_i\not=\tau^i x'_i$ and thus that the sectional path from $\tau^my$
to $x$ provided by the above argument is also on a different sequence
of vertices. Thus, applying the above argument to every
sectional path of length $m$ from $x$ to $y$, we obtain an injection from
the set of sectional paths of length $m$ from $x$ to $y$ to the set of
sectional paths of length $m$ from $\tau^my$ to $x$.

A similar argument shows that whenever there is 
a sectional path 
\[
\tau^my=y_0\to y_1\to\dots\to y_m=x
\]
in $\Gamma$ with multiplicities $l_1,l_2,\ldots ,l_m$,
then $\tau^{i-m}y_i$ is defined for all $i$ and
there is a sectional path 
\[
x\to\tau^{-1}y_{m-1}\to\dots\to\tau^{m-1}y_1\to\tau^my=y_0
\]
in $\Gamma$ with multiplicities $l_m,l_{m-1},\ldots ,l_1$
and as above we obtain an injection from the set of sectional paths of
length $m$ from $\tau^my$ to $x$ to the set of sectional paths of length
$m$ from $x$ to $y$.

Since $\Gamma$ is locally finite, the number of sectional paths of fixed
length between two vertices is finite. It follows that the number of sectional
paths of length $m$ from $x$ to $y$ is the same as the number of sectional
paths of length $m$ from $\tau^my$ to $x$.
Hence $(\Gamma^m,\tau^m)$ is a translation quiver.
\end{proof}

We remark that the square of the translation quiver below, which does not
satisfy the additional assumption of the theorem, is not a translation quiver:
$$
\xymatrix@-4mm{
&& {\circ} \ar[rd] & & \\
&&& {\circ} \ar[rd] & \\
\circ \ar@{.}[rr] && \circ \ar[ur] \ar@{.}[rr] && \circ
}
$$

\begin{cor}
(1) Let $(\Gamma,\tau)$ be a hereditary translation quiver. Then
$(\Gamma^m,\tau^m)$ is a translation quiver. \\
(2) Let $(\Gamma,\tau)$ be a stable translation quiver. Then
$(\Gamma^m,\tau^m)$ is a stable translation quiver.
\end{cor}

\begin{proof}
Part (1) is immediate from Theorem~\ref{mpower} and the definition of
a hereditary translation quiver. For (2), note that
if $(\Gamma,\tau)$ is stable, no vertex is projective, so $(\Gamma^m,\tau^m)$
is a translation quiver by Theorem~\ref{mpower}.
Since $\tau$ is defined on all vertices of $\Gamma$, so is $\tau^m$.
\end{proof}

We remark that the $m$th power of a hereditary translation quiver need not
be hereditary: there can be non-projective vertices $z$ without any vertex
$z'$ such that $z'\to z$. For example, consider the hereditary translation
quiver below. It is clear that its square in the above sense has no arrows,
but does have non-projective vertices.

$$
\xymatrix@-5mm{
 & \ar[rd] \ar@{.}[rr] && \ar[rd] \ar@{.}[rr] && \ar[rd] && {\cdots} \\
\ar[ur] \ar@{.}[rr] && \ar[ur] \ar@{.}[rr] && \ar[ur] \ar@{.}[rr] &&& {\cdots}
}
$$

However, we do have the following:

\begin{prop}
Let $(\Gamma,\tau)$ be a translation quiver such that for any
arrow $x\to y$ in $\Gamma$, $x$ is projective whenever
$y$ is projective. Then the translation quiver $(\Gamma^m,\tau^m)$ has
the same property.
\end{prop}

\begin{proof}
We know by Theorem~\ref{mpower} that $(\Gamma^m,\tau^m)$ is a translation
quiver. Suppose that $$x_0=x\to x_1\to \cdots \to x_m=y$$
is a sectional path in $\Gamma$ and that $\tau^m x$ is defined.
We show by induction on $m$ that $\tau^m y$ is defined.
It is clear that this holds for $m=1$, since $\Gamma$ is hereditary.
Since $x$ is not projective, $x_1,x_2,\ldots ,x_m$ are not projective.
Since there are arrows $x_{i-1}\to x_{i}$ for $i=1,2,\ldots ,m$,
there are arrows $\tau x_{i}\to x_{i-1}$ and therefore arrows
$\tau x_{i-1}\to \tau x_i$ for $i=1,2,\ldots ,m$. It is clear
that any path
$$\tau x_0=\tau x\to \tau x_1\to \cdots \to \tau x_m=y$$
is sectional; also $\tau^{m-1}(\tau x)$ is defined. Hence, by the inductive
hypothesis, $\tau^{m-1}(\tau x_m)=\tau^m x_m$
is defined and we are done.
\end{proof}

\section{The $m$-cluster category in terms of $m$th powers}
We consider the construction of Section~\ref{se:mpower}
in the case where $\Gamma$ is the quiver 
given by the diagonals of an $N$-gon $\Pi$, 
i.e. $\Gamma=\Gamma^1_{A_{N-3}}$ as in Section~\ref{se:quiver}.
Here, we fix $m=1$, i.e. 
the vertices of the quiver are the usual diagonals 
of $\Pi$ and there is an arrow from $D$ to $D'$ 
if $D,D'$ have a common endpoint $i$ so that
$D,D'$ together with the arc from $j$ to $j'$ between the 
other endpoints form a triangle and $D$ is rotated to
$D'$ by a clockwise rotation about $i$.
We will call this rotation $\rho_i$. 
Furthermore, we have introduced an automorphism 
$\tau_1$ of $\Gamma$: $\tau_1$ 
sends $D$ to $D'$ if $D$ can be rotated  
to $D'$ by an anticlockwise rotation about the centre of 
the polygon through $\frac{2\pi}{N}$. Then 
$\Gamma=\Gamma^1_{A_{N-3}}$ is a stable 
translation quiver (cf. Proposition~\ref{prop:stable}). 

The geometric interpretation of a sectional path 
of length $m$ from $D$ to $D'$ is 
given by the map $\rho_i^m$: 
$\rho_i^m$ sends 
the diagonal $D$ to $D'$ if $D$, $D'$ have 
a common endpoint $i$ and form an $m+2$-gon 
together with the arc between the other 
endpoints $j,j'$ and if $D$ can be rotated 
to $D'$ with a clockwise rotation about the common 
endpoint.  

Furthermore, the $m$-th power $\tau_1^m$ of the translation 
$\tau_1$ corresponds 
to a anticlockwise rotation through $\frac{2m\pi}{N}$ about the 
centre of the polygon. From that one obtains: 

\begin{prop}\label{prop:m-quiver}
1)
The quiver $(\Gamma^1_{A_{N-3}})^m$ contains a 
translation quiver of $m$-diagonals if and only 
if $N=nm+2$ for some $n$. 

2) $\Gamma^m_{A_{n-1}}$ is a 
connected component of $(\Gamma^1_{A_{nm-1}})^m$. 
\end{prop}

\begin{proof}
1) Note that if $N\neq nm+2$ (for some $n$) 
then $\Gamma^m_{A_{N-3}}$ contains no $m$-diagonals. 

So assume that $N=nm+2$ for some $n$. 
Let $\Gamma:=\Gamma^1_{A_{N-3}}=\Gamma^1_{A_{nm-1}}$. 
We have to show that $\Gamma^m$ 
contains $Q:=\Gamma^m_{A_{n-1}}$. 
Recall that the vertices of the quiver 
$\Gamma^m$ are the diagonals of an $nm+2$-gon 
and that 
$Q$ is the quiver 
whose vertices are the $m$-diagonals of an $nm+2$-gon.
So the vertices of $Q$ are vertices of $\Gamma^m$. 

We claim that the arrows between those vertices are the same 
for $Q$ and for $\Gamma^m$. In other words, we claim that there 
is a sectional path of length $m$ between $D$ and $D'$ if 
and only if $D$ can be rotated clockwise to $D'$ about a common 
endpoint and $D$ and $D'$ together with an arc joining the other 
endpoints bound an $(m+2)$-gon. 

Let $D\to D'$ be an arrow in $\Gamma^m$, where
$D$ is the diagonal $(i,j)$ from $i$ to $j$.
Without loss of generality, let $i<j$. 
The arrow $D\to D'$ in $\Gamma^m$ corresponds to a 
sectional path of length $m$ in $\Gamma$, 
$D\to D_1\to\dots\to D_{m-1}\to D_m=D'$. 
We describe such sectional paths. The first arrow 
is either $D=(i,j)\to (i,j+1)$
or $(i,j)\to(i-1,j)$, i.e. 
$D_1=(i,j+1)$ or $D_1=(i-1,j)$ (vertices 
taken $\mod N$). 
In the first case, one then gets an arrow 
$D_1=(i,j+1)\to(i+1,j+1)$ or $D_1\to(i,j+2)$. 
Now $\tau(i+1,j+1)=(i,j)$ and since the 
path is sectional, we get that $D_2$ can 
only be the diagonal $(i,j+2)$. 
$$
{\small
\xymatrix@-8mm{
%\xymatrix{
 & & & & (i,j+2)&\\
 & & (i,j+1)%\ar[rrd]
\ar[rru] & \\
(i,j)\ar[rru]\ar[rrd]\ar@{.}[rrrr] & & & & (i+1,j+1)\\
 & & (i-1,j)%\ar[rru]
\ar[rrd] & \\
 & & & & (i-2,j)
}}
$$
Repeatedly using the above argument, we see that the sectional path has to 
be of the form 
\[
D=D_0=(i,j)\to (i,j+1)\to (i,j+2)\to\dots\to(i,j+m)=D_m=D' 
\]
where all vertices are taken $\mod N$.

Similarly, if $D_1=(i-1,j)$ then $D_2=(i-2,j)$ 
and so on, $D_m=(i-m,j)$ ($\mod N$). 

In particular, in the first case, the arrow $D\to D'$ corresponds to a rotation 
$\rho_i^m$ about the common
endpoint $i$ of $D,D'$. In the second case, the arrow 
$D\to D'$ corresponds to $\rho_j^m$. 
In each case $D, D'$ and an arc between them bound an $(m+2)$-gon, so
there is an arrow from $D$ to $D'$ in $Q$.

Since it is clear that every arrow in $Q$ arises in this way, we
see that the arrows between the vertices of $Q$ and of the corresponding
subquiver of $\Gamma^m$ are the same. 

2) We know by Proposition~\ref{prop:connected} that 
$Q=\Gamma^m_{A_{n-1}}$ is a connected stable translation quiver. 
If there is an arrow $D\to D'$ in $\Gamma^m$ where $D$ is an 
$m$-diagonal then $D'$ is an $m$-diagonal. Similarly, 
$\tau_1^m(D)$ is also an $m$-diagonal.
\end{proof}

\begin{theorem}
The $m$-cluster category $\mathcal{C}^m_{A_{n-1}}$ 
is a full subcategory of the additive category generated 
by the mesh category of $(\Gamma^1_{A_{nm-1}})^m$. 
\end{theorem}
\begin{proof}
This is a consequence of Proposition~\ref{prop:m-quiver} 
and Theorem~\ref{thm:m-cluster}
\end{proof}

\begin{remark}
Even if $\Gamma$ is a connected quiver, 
$\Gamma^m$ need not be connected. As an 
example we consider the quiver 
$\Gamma=\Gamma^1_{A_5}$ and its second 
power $(\Gamma^1_{A_5})^2$ pictured in 
Figures~\ref{fig:Gamma1-1} and~\ref{fig:Gamma1-2}.
The connected components of 
$(\Gamma^1_{A_5})^2$ are $\Gamma^2_{A_2}$ 
and two copies of a translation quiver whose mesh 
category is equivalent to ${\rm ind} D^b(A_3)/[1]$
(where $D^b(A_3)$ denotes the derived category of a
Dynkin quiver of type $A_3$).
We thus obtain a geometric construction of a quotient
of $D^b(A_3)$ which is not an $m$-cluster category.

\begin{figure}[h] 
\begin{center}
\includegraphics{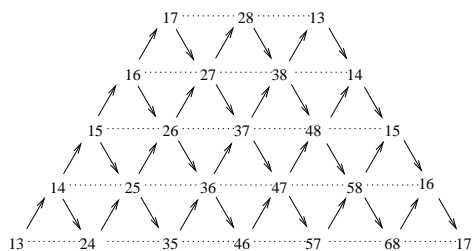}
\caption{The quiver $\Gamma^1_{A_5}$} \label{fig:Gamma1-1}
\end{center}
\end{figure} 

\begin{figure}[h] 
\begin{center}
\includegraphics{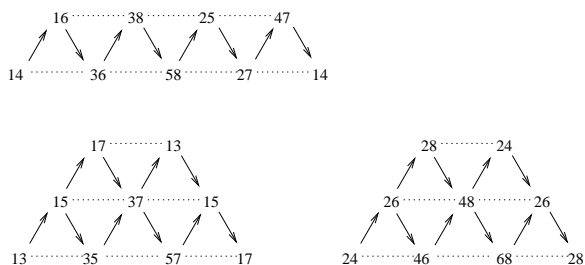}
\caption{The three components of $(\Gamma^1_{A_5})^2$} \label{fig:Gamma1-2}
\end{center}
\end{figure} 
\end{remark}

%%%%%%%%%%%%%%%%%%%%%%%%%%%%%%%%%%%%%%%%%%
%


\begin{thebibliography}{99}

\bibitem[BMRRT]{bmrrt}
A. B. Buan, R. J. Marsh, M. Reineke, I. Reiten and G. Todorov. 
\textit{Tilting theory and cluster combinatorics.}
Preprint arxiv:math.RT/0402054, February 2004, to appear in Adv. Math.  

\bibitem[CCS1]{ccs1} P. Caldero, F. Chapoton and R. Schiffler, 
\textit{Quivers with Relations arising from Clusters 
($A_n$ case).} Trans. Amer. Math. Soc.  358  (2006),  
no. 3, 1347--1364. 

\bibitem[CCS2]{ccs2}
P. Caldero, F. Chapoton and R. Schiffler,
\textit{Quivers with relations and cluster tilted algebras.}
Preprint arxiv:math.RT/0411238, November 2004, to appear in J. Alg.
Rep. Theory.

\bibitem[FR]{fominreading}
S. Fomin and N. Reading.
\emph{Generalized cluster complexes and Coxeter combinatorics.}
International Mathematics Research Notices \textbf{2005}:44, 2709-2757. 

\bibitem[FZ1]{fominzelevinsky1}
S. Fomin and A. Zelevinsky,
\textit{Cluster algebras I: Foundations.}
J. Amer. Math. Soc. \textbf{15} (2002), no. 2, 497--529. 

\bibitem[FZ2]{fominzelevinsky2}
S. Fomin and A. Zelevinsky,
\textit{Y-systems and generalized associahedra.}
Annals of Mathematics.  \textbf{158} (2003), No. 3,  977--1018. 

\bibitem[Hap]{happel}
D. Happel.
\textit{Triangulated categories in the representation theory of finite
dimensional algebras.}
LMS Lecture Note Series 119, Cambridge University Press, 1988.

\bibitem[Kel]{keller}
B. Keller.
\textit{On triangulated orbit categories.}
Documenta Math. \textbf{10} (2005), 551-581.

\bibitem[KR]{kellerreiten}
B. Keller and I. Reiten.
\textit{Cluster-tilted algebras are Gorenstein and stably Calabi-Yau.}
Preprint arxiv:math.RT/0512471, December 2005.

\bibitem[Rie]{riedtmann}
C. Riedtmann.
\textit{Algebren, Darstellungsk\"{o}cher, \"{U}berlagerungen und zur\"{u}ck.}
Comment. Math. Helv. \textbf{55} (1980), no. 2, 199--224.

\bibitem[Rin]{ringel}
C.M. Ringel, 
\textit{Tame Algebras and Integral Quadratic 
Forms.} Lecture Notes in Mathematics 
\textbf{1099} (1984), Springer, Berlin.

\bibitem[Tho]{thomas}
H. Thomas.
\textit{Defining an $m$-cluster category.}
Preprint, 2005.

\bibitem[Tza]{tzanaki}
E. Tzanaki.
\textit{Polygon dissections and some generalizations of cluster complexes.}
Preprint arxiv:math.CO/0501100, January 2005.

\bibitem[Zhu]{zhu}
B. Zhu.
\textit{Generalized cluster complexes via quiver representations.}
Preprint, July 2006.
\end{thebibliography}
\end{document}